\documentclass[11pt]{amsart}
 \usepackage{amscd,amsmath,amsthm,amssymb}
 \usepackage{pstcol,pst-plot,pst-3d}
 \psset{unit=0.7cm,linewidth=0.8pt,arrowsize=2.5pt 4}

\newpsstyle{fatline}{linewidth=1.5pt}
\newpsstyle{fyp}{fillstyle=solid,fillcolor=verylight}
\definecolor{verylight}{gray}{0.97}
\definecolor{light}{gray}{0.93}
\definecolor{medium}{gray}{0.82}
\definecolor{dark}{gray}{0.72}
 \def\opn#1#2{\def#1{\operatorname{#2}}} 
 \opn\chara{char} \opn\length{\ell} \opn\pd{pd} \opn\rk{rk}
 \opn\projdim{proj\,dim} \opn\injdim{inj\,dim} \opn\rank{rank}
 \opn\depth{depth} \opn\grade{grade} \opn\height{height}
 \opn\embdim{emb\,dim} \opn\codim{codim}
 
 \opn\Tr{Tr} \opn\bigrank{big\,rank}
 \opn\superheight{superheight}\opn\lcm{lcm}
 \opn\trdeg{tr\,deg}
 \opn\reg{reg} \opn\lreg{lreg} \opn\ini{in} \opn\lpd{lpd}
 \opn\size{size} \opn\sdepth{sdepth}
 \opn\link{link}\opn\fdepth{fdepth}\opn\lex{lex}
 \opn\div{div} \opn\Div{Div} \opn\cl{cl} \opn\Cl{Cl}
 \opn\Spec{Spec} \opn\Supp{Supp} \opn\supp{supp} \opn\Sing{Sing}
 \opn\Ass{Ass} \opn\Min{Min}\opn\Mon{Mon}
 \opn\Ann{Ann} \opn\Rad{Rad} \opn\Soc{Soc}
 \opn\Im{Im} \opn\Ker{Ker} \opn\Coker{Coker} \opn\Am{Am}
 \opn\Hom{Hom} \opn\Tor{Tor} \opn\Ext{Ext} \opn\End{End}
 \opn\Aut{Aut} \opn\id{id}
 
 \opn\nat{nat}
 \opn\pff{pf}
 \opn\Pf{Pf} \opn\GL{GL} \opn\SL{SL} \opn\mod{mod} \opn\ord{ord}
 \opn\Gin{Gin} \opn\Hilb{Hilb}\opn\sort{sort}
 \opn\aff{aff} \opn\con{conv} \opn\relint{relint} \opn\st{st}
 \opn\lk{lk} \opn\cn{cn} \opn\core{core} \opn\vol{vol}
 \opn\link{link} \opn\star{star}\opn\lex{lex}\opn\set{set}
 \opn\gr{gr}
 
 \def\pot#1#2{#1[\kern-0.28ex[#2]\kern-0.28ex]}

 \opn\dirlim{\underrightarrow{\lim}}
 \opn\inivlim{\underleftarrow{\lim}}
 \let\union=\cup
 \let\sect=\cap

 \def\Implies{\ifmmode\Longrightarrow \else
         \unskip${}\Longrightarrow{}$\ignorespaces\fi}
 \def\implies{\ifmmode\Rightarrow \else
         \unskip${}\Rightarrow{}$\ignorespaces\fi}
 \def\iff{\ifmmode\Longleftrightarrow \else
         \unskip${}\Longleftrightarrow{}$\ignorespaces\fi}

 \let\:=\colon
 \newtheorem{Theorem}{Theorem}[section]
 \newtheorem{Lemma}[Theorem]{Lemma}

 \let\epsilon\varepsilon
 \let\kappa=\varkappa
 \def\qed{\ifhmode\textqed\fi
       \ifmmode\ifinner\quad\qedsymbol\else\dispqed\fi\fi}
 \def\textqed{\unskip\nobreak\penalty50
        \hskip2em\hbox{}\nobreak\hfil\qedsymbol
        \parfillskip=0pt \finalhyphendemerits=0}
 \def\dispqed{\rlap{\qquad\qedsymbol}}
 
 \opn\dis{dis}
 \def\pnt{{\raise0.5mm\hbox{\large\bf.}}}
 
 \opn\Lex{Lex}

\begin{document}

 \title {On the stable set of associated prime ideals of a monomial ideal}

 \author {Shamila Bayati, J\"urgen Herzog and Giancarlo Rinaldo}

\address{Shamila Bayati, Faculty of Mathematics and Computer Science,
Amirkabir University of Technology (Tehran Polytechnic), 424 Hafez Ave., Tehran
15914, Iran.
}\email{shamilabayati@gmail.com}

\address{J\"urgen Herzog, Fachbereich Mathematik, Universit\"at Duisburg-Essen, Campus Essen, 45117
Essen, Germany} \email{juergen.herzog@uni-essen.de}

\address{Giancarlo Rinaldo, Dipartimento di Matematica\\
Universit\`a di Messina\\
Contrada Papardo, salita Sperone, 31\\
98166 Messina, Italy}
\email{giancarlo.rinaldo@tiscali.it}

 \begin{abstract}
 It is shown that any set of nonzero monomial prime ideals can be realized as the stable set of associated prime ideals of a monomial ideal. Moreover, an algorithm is given to compute the stable set of associated prime ideals of a monomial ideal.
 \end{abstract}

\thanks{
This paper was written while S.\ Bayati and G.\ Rinaldo  were visiting Universit\"at Duisburg-Essen, Campus Essen.
The first author was supported by the Ministery of Science of Iran, the third author by INDAM GNSAGA}

 \maketitle

 \section*{Introduction}

It is known by a result of Brodmann \cite{B} that in any Noetherian ring the set of associated prime ideals  $\Ass(I^s)$ for the powers of an ideal $I$ stabilizes for
$s\gg 0$. In other words, there exists an integer $s_0$ such that $\Ass(I^s)=\Ass(I^{s+1})$ for all $s\geq s_0$. This stable set of associated prime ideals is denoted by $\Ass^\infty(I)$. In recent years there have been several publications \cite{CMS,FHV,HRV,MMV}, mostly in combinatorial contexts,  to describe $\Ass^\infty(I)$.

The main purpose of this note is to show that for any set of nonzero monomial  prime ideals,  there exists a monomial  ideal for which  this given  set is precisely the set  of stable associated prime ideals. This result is given in Theorem~\ref{ourlasthope}. Describing  the possible stable sets of associated prime ideals for  squarefree monomial ideals remains an open problem.

In Section 2 we give an algorithm that determines  $\Ass^\infty(I)$ for any monomial ideal. The routine written in {\em Macaulay 2} can be found in \cite{BHR}.

\section{The set $\Ass^\infty(I)$ of a monomial ideal}

For the proof of the main result of this section we first need the following general fact about associated prime ideals.

\begin{Lemma}
\label{shamila139}
Let $R$ be a Noetherian ring,  $I\subset R$ an ideal and $P\subset R$ a prime ideal such that $P\not\subseteq  P'$ for all $P'\in \Ass(I)$. Let $k$ be an  integer such that $I\not\subseteq P^{(k)}$, where $ P^{(k)}$ denotes  the $k$th symbolic power of $P$.   Then $\Ass(I\sect P^{(k)})=\Ass(I)\union\{P\}$.
\end{Lemma}
\begin{proof}
Let $I=Q_1\cap Q_2\cap \ldots \cap Q_m$ be a irredundant primary decomposition of $I$ with $\sqrt{Q_i}=P_i$ and $P_i\neq P_j$ for $i\neq j$. We show that $I\sect P^{(k)}=Q_1\cap Q_2\cap \ldots \cap Q_m\cap P^{(k)}$ is a irredundant  primary decomposition of $I\sect P^{(k)}$. Considering the fact $I\not\subseteq P^{(k)}$, it is enough to show that
\begin{eqnarray}
\label{goodenglish}
(\bigcap_{j\neq i}Q_j)\cap P^{(k)}\not\subseteq Q_i
\end{eqnarray}
for all $i$. Since  $I=Q_1\cap Q_2\cap \ldots \cap Q_m$ is a irredundant primary decomposition of $I$, there exists an element $a\in  (\bigcap_{j\neq i}Q_j)\setminus Q_i$. There also exists an element $b\in P\setminus P_i$ by our assumption. So $ab^k\in (\bigcap_{j\neq i}Q_j)\cap P^{(k)}$, but  $ab^k\not\in Q_i$. Indeed, if $ab^k\in Q_i$ it follows that $a\in Q_i$ or $b^k\in \sqrt{Q_i}=P_i$, a contradiction. This implies (\ref{goodenglish}).
\end{proof}

Now we consider a more specific situation. Let $K$ be a field, $S=K[x_1,\ldots,x_n]$ the polynomial ring in the variables $x_1,\ldots,x_n$ and $I\subset S$ a monomial ideal. As usual we denote by $G(I)$ the unique minimal set of monomial generators of $I$.

The associated prime ideals of $I$ are all monomial prime ideals, that is, ideals of the form $P_F=(\{x_i\:\ i\in F\})$ where $F\subset [n]$.  The main result of this section is the following

 \begin{Theorem}
 \label{ourlasthope}
 Let $P_1,\ldots,P_m\subset S$ be an arbitrary collection of nonzero monomial prime ideals. Then there exists a monomial ideal $I\subset S$ such that
 \[
 \Ass^\infty(I)=\{P_1,\ldots,P_m\}.
 \]
 \end{Theorem}

 \begin{proof}  We may assume that $P_i\neq P_j$ for $i\neq j$, and that $|G(P_i)|\leq |G(P_j)|$ for $i<j$. For $r=1,\ldots,m$, we inductively  define  the  monomial ideals $J_r$ with $\Ass(J_r)=\{P_1,\ldots.P_r\}$ as follows: $J_1=P_1$. Suppose that for some $r<m$ the ideal $J_r$ is already defined. Then we set
 \[
 J_{r+1}=J_r\sect P_{r+1}^{k_{r+1}},
 \]
 where $k_{r+1}>\min\{\deg u\:\; u\in G(J_r)\}$. By Lemma~\ref{shamila139} we have that
 \[
 \Ass(J_{r+1})=\{P_1,\ldots,P_{r+1}\}.
 \]
 Indeed, since $|G(P_i)|\leq |G(P_j)|$ for $i<j$, it follows that $P_{r+1}\not\subseteq P_i$ for $i\leq r$. Moreover, $P_{r+1}^{(k_{r+1})}=P_{r+1}^{k_{r+1}}\not\supseteq J_r$,  by the choice of $k_{r+1}$.

 Let $t$ be any positive integer. Then
 $tk_{r+1}$ is bigger than the minimum degree  of $P_1^{tk_1}\sect P_2^{tk_2}\sect \ldots\sect  P_r^{tk_r}$ for $r=1,\ldots,m-1$, because $J_r^t\subset P_1^{tk_1}\sect P_2^{tk_2}\sect \ldots\sect  P_r^{tk_r}$. Thus as before we see that
 \begin{eqnarray}
 \label{t}
 \Ass( P_1^{tk_1}\sect P_2^{tk_2}\sect \ldots\sect  P_m^{tk_m})=\{P_1,\ldots,P_m\} \quad \text{for all} \quad t.
 \end{eqnarray}
By \cite[Corollary 2.2]{HHT},  there exists an integer $d$ such that for
\[
I=P_1^{dk_1}\sect P_2^{dk_2}\sect \ldots\sect  P_m^{dk_m}
\]
we have
\[
I^s=P_1^{sdk_1}\sect P_2^{dk_2}\sect \ldots\sect  P_m^{sdk_m}\quad \text{for all} \quad s\geq 1.
\]
Thus (\ref{t}) implies that $\Ass(I)=\Ass(I^s)=\{P_1,\ldots,P_m\}$. This yields the desired conclusion.
 \end{proof}

The monomial ideal $I$ that we constructed in Theorem~\ref{ourlasthope} has the property that $\Ass(I)=\Ass^\infty(I)$. In general, for any  ideal $I$ one has $\Min(I)=\Min(I^s)$ for all $s$. Hence as strengthening of Theorem~\ref{ourlasthope},  one could ask the following question: suppose we are given two sets $A=\{P_1,\ldots,P_\ell\}$ and $B=\{P_1',\ldots, P_m'\}$  of monomial prime ideals such that the minimal elements of these sets with respect to inclusion are the same. For which such sets  does exist a monomial  ideal $I$ such that $\Ass(I)=A$ and $\Ass^\infty(I)=B$?   For example, there is no monomial ideal $I$ with $\Ass(I)=\{(x_1),(x_2)\}$ and $\Ass^\infty(I)=\{(x_1),(x_2),(x_1,x_2)\}$.

\medskip
It would be also interesting to understand the possible sets of prime ideals for $\Ass^\infty(I)$ when $I$ is a {\em squarefree} monomial ideal. Consider for example the set of monomial prime ideals $\{(x_1),(x_2),(x_1,x_2)\}$. This set cannot be $\Ass^\infty(I)$ for any squarefree monomial ideal $I$. Thus there is no analogue of Theorem~\ref{ourlasthope} for squarefree monomial ideals.

 \section{An algorithm to compute $\Ass^\infty(I)$}

 In this section we describe an algorithm that we implemented in  {\em Macaulay 2} (see \cite{BHR}) to compute $\Ass^\infty(I)$  for a monomial ideal $I$. Since there are only finitely many monomial prime ideals in $S$, it is enough to have an algorithm to determine  whether  a monomial prime ideal $P\subset S$  belongs to $\Ass^\infty(I)$.

 The input of our algorithm  is a monomial ideal $I$ given by its set of monomial generators $G(I)$,   and a monomial prime ideal $P=(x_{i_1},\ldots,x_{i_k})$.

 \begin{enumerate}
 \item[ Step 1:] Define  the ideal $J\subset S'=K[x_{i_1},\ldots,x_{i_k}]$ generated by the monomials $u^*$ with $u\in G(I)$, where $u^*$ is obtained from $u$ by the substitution $x_j\mapsto 1$ for $x_j\not \in P$.

\item[  Step 2:] If $J=S'$, then $P\not\in \Ass^\infty(I)$ (actually $I\not\subseteq P$). Else go to Step 3.

\item[  Step 3:] Form the Rees algebra ${\mathcal R}(J)$ and compute the Koszul homology $$H= H_{k-1}(x_{i_1},\ldots,x_{i_k};{\mathcal R}(J)).$$

\item[  Step 4:] Compute the Krull dimension  of $H$. If $\dim H>0$, then  $P\in \Ass^\infty(I)$; else $P\not\in \Ass^\infty(I)$.
 \end{enumerate}

\medskip
\noindent
Indeed, to justify this algorithm notice that  $P\in \Ass^\infty(I)$ if and only if $PS_x\in \Ass(J^sS_x)$ for all  $s\gg 0$, where $x=\prod_{ x_i\not\in P}x_i$. This is equivalent to say that $\depth S'/J^s=0$ for all $s\gg 0$. Next we observe that the $j$-th Koszul homology  $H_j(x_{i_1},\ldots,x_{i_k};\mathcal(J))$ is a graded ${\mathcal R}(J)$-algebra with
\[
H_j(x_{i_1},\ldots,x_{i_k};{\mathcal R}(J))_s=H_j(x_{i_1},\ldots,x_{i_k};J^s) \quad \text{for all} \quad s.
\]
Hence we have $\depth S'/J^s=0$ if and only if  $H_{k-1}(x_{i_1},\ldots,x_{i_k};{\mathcal R}(J))_s \neq 0$. Therefore   $\depth S'/J^s=0$ for all $s\gg 0$ if and only if the  finitely generated ${\mathcal R}(J)$-module $H=H_{k-1}(x_{i_1},\ldots,x_{i_k};{\mathcal R}(J))$ has infinitely many nonzero components, and this is the case if and only the Krull dimension of $H$ is positive.

\medskip
To demonstrate this algorithm we consider the following computation in {\em Macaulay 2}. We want to check whether  $P=(a,b,c,e,f)\in \Ass^\infty(I)$ for $I=(a^3b^2e,bc^3d,b^4de^2f, ab^2cf^3)$.

\medskip
\begin{verbatim}
S=QQ[a..f];
J=monomialIdeal(a^3*b^2*e,b*c^3,b^4*e^2*f, a*b^2*c*f^3);
R=reesAlgebra(J);
P={a,b,c,e,f};
phi=matrix{P};
C=(koszul phi)**R;
dim(HH_4 C)>0
\end{verbatim}

\medskip

\noindent
In this case the output is: {\tt true}.  This  means that $P\in \Ass^\infty(I)$. Other examples and the source file for a routine can be found in \cite{BHR} where $\Ass^\infty(I)$ is computed for monomial ideals.

{}

\begin{thebibliography}{}
\bibitem{B} M. Brodmann, Asymptotic stability of ${\Ass}(M/I^{n}M)$,
Proc. Amer. Math. Soc. 74 (1979), no. 1, 16--18.

\bibitem{BHR} S.\ Bayati, J.\ Herzog and G.\ Rinaldo, A routine to compute the stable set of associated prime ideals of a monomial ideal,
{\tt Available at http://ww2.unime.it/algebra/rinaldo/stableset/}, 2011.

\bibitem{CMS} J. Chen, S. Morey, A. Sung, The stable set of associated primes of the ideal of a graph,
Rocky Mountain J. Math. 32 (2002), no. 1, 71--89.

\bibitem{FHV} C. A. Francisco, H. T\`{a}i H\`{a} A. Van Tuyl,  Colorings of hypergraphs, perfect graphs, and associated primes of powers of monomial ideals, J. Algebra 331 (2011), 224--242.

\bibitem{GS} D. R. Grayson, M. E. Stillman, Macaulay 2, a software system for research in algebraic geometry.
http://www.math.uiuc.edu/Macaulay2/.

\bibitem{HHBook} J. Herzog, T. Hibi,  Monomial ideals. Graduate Texts in Mathematics, 260. Springer-Verlag London, Ltd., London, 2011.


\bibitem{HHT} J. Herzog, T. Hibi, N. V. Trung,  Symbolic powers of monomial ideals and vertex cover algebras, Adv. Math. 210 (2007), no. 1, 304--322.

\bibitem{HRV} J. Herzog, A. Rauf, M. Vladoiu, The stable set of associated prime ideals of a polymatroidal ideal, Preprint, 2011,
arXiv:1109.5834v2 [math.AC].

\bibitem{MMV} J. Martinez-Bernal, S. Morey, R. H. Villarreal, Associated primes of powers of edge ideals, Preprint, 2011,
arXiv:1103.0992v3 [math.AC].




\end{thebibliography}
 \end{document}